\documentclass[oneside,english]{amsart}
\usepackage[latin9]{inputenc}
\usepackage{geometry}
\usepackage{babel}
\usepackage{varwidth}
\usepackage{amsbsy}
\usepackage{amstext}
\usepackage{amsthm}
\usepackage{amssymb}
\usepackage{setspace}
\usepackage{microtype}
\usepackage[bookmarks=false,
 breaklinks=false,pdfborder={0 0 0},pdfborderstyle={},backref=false,colorlinks=false]
 {hyperref}

\makeatletter

\providecommand{\tabularnewline}{\\}
\newenvironment{cellvarwidth}[1][t]
    {\begin{varwidth}[#1]{\linewidth}}
    {\@finalstrut\@arstrutbox\end{varwidth}}

\numberwithin{equation}{section}
\numberwithin{figure}{section}
\theoremstyle{plain}
\newtheorem{thm}{\protect\theoremname}
\theoremstyle{remark}
\newtheorem{rem}[thm]{\protect\remarkname}
\theoremstyle{plain}
\newtheorem{lem}[thm]{\protect\lemmaname}

\@ifundefined{date}{}{\date{}}
\usepackage{calrsfs}
\let\cal=\mathcal

\def\noal#1{\noalign{\hbox{#1}}}

\makeatother

\providecommand{\lemmaname}{Lemma}
\providecommand{\remarkname}{Remark}
\providecommand{\theoremname}{Theorem}

\begin{document}
\title{Circularity in Finite Fields and Solutions of the Equations $\boldsymbol{x^{m}+y^{m}-z^{m}=1}$}
\author{Wen-Fong Ke}
\address{Department of Mathematics, National Cheng Kung University, Tainan
701, Taiwan}
\email{wfke@mail.ncku.edu.tw}
\author{Hubert Kiechle}
\address{Mathematisches Seminar, Universit\"at Hamburg, Bundesstr. 55, D-20146
Hamburg, Germany}
\email{hubert.kiechle@uni-hamburg.de}
\begin{abstract}
An explicit formula for the number of solutions of the equation in
the title is given when a certain condition, depending only on the
exponent and the characteristic of the field, holds. This formula
improves the one given by the authors in an earlier paper.
\end{abstract}

\maketitle

\section{Introduction}

\global\long\def\vp{\varphi}%
\global\long\def\bk{\boldsymbol{k}}%
\global\long\def\cN{\mathcal{N}}%
\global\long\def\cC{\mathcal{C}}%
\global\long\def\bZ{\mathbb{Z}}%
\global\long\def\bC{\mathbb{C}}%
\global\long\def\bQ{\mathbb{Q}}%
\global\long\def\cPk{{\cal P}_{k}}%
Let $p,q,r,k\in\mathbb{N}$ with $p$ prime, $k\geq3$ and $q=p^{r}$
such that $k\mid(q-1)$. Put $m=(q-1)/k$. In \cite{KK95} and \cite{K94a}
an explicit formula is given for the number of solutions of the equation
\begin{equation}
x^{m}+y^{m}-z^{m}=1\label{eq:maxp}
\end{equation}
in the Galois field $F=\text{GF}(q)$ under the following hypothesis:
the unique multiplicative subgroup $\Phi$ of order $k$ in $F$ satisfies
$|\Phi a\cap(\Phi b+c)|\leq2$ for all $a,b,c\in F^{*}=F\setminus\{0\}$.
In this case, the pair $(F,\Phi)$ (known as a \textit{Ferrero pair},
see \cite{Clay92}) as well as the pair $(q,k)$, is said to be \textit{circular.}

Equations of similar type over finite fields have been studied for
a long time; see \cite{Joly73} for a comprehensive exposition. In
general, finding explicit formulas for such equations is difficult,
and has only been achieved in special cases. Namely, for small $m$,
eg., \cite{Dickson36,Lehmer91}, for small $k$, eg., \cite{Lehmer91,Sun89},
and for some other restrictive conditions as in \cite{AnuradhaK99,Faircloth52,HuaVandiver49b,Wolfmann92}.
Circularity restricts the size of $k$ and implies a lower bound for
$m$. Indeed, by \cite[(4.1)]{KK95} we have that if $(q,k)$ is circular,
then
\[
m\ge\frac{\sqrt{4q-7}-3}{2}.
\]
Thus we attack from the top.

It turns out that circularity of $(q,k)$ dependents only on $p$
and $k$ (\cite[Theorem~5]{Modisett89}, see also \cite[(5.21)]{Clay92}).
That is, $(q,k)$ is circular if and only if $(p^{r'},k)$ is circular
for some $r'\in\mathbb{N}$ such that $k\mid(p^{r'}-1)$. Thus, we
say that the pair $(p,k)$ is circular if $(p^{r'},k)$ is circular
for some $r'\in\mathbb{N}$.

Now, the formula for the number of solutions of (\ref{eq:maxp}) given
in \cite{KK95} and \cite{K94a} was actually under another constraint
when $k$ is odd. Namely, when $k$ is odd, it was conveniently required
that $(p,2k)$ is also circular. This excludes the case of characteristic~$2$
for the field.

In this article, we investigate the circularity of $(p,k)$ more closely
to remove the extra condition imposed on the formula. Effectively,
we will have the following theorem.
\begin{thm}
\label{Thm1_KK95} Let $(q,k)$ be circular, also, let $N'$ be the
number of solutions with $xyz\not=0$.
\begin{enumerate}
\item If $k$ is even, then
\[
N=\begin{cases}
3(k-1)m^{3}+6m^{2}+3m, & \text{if \ensuremath{6\mid k};}\\
3(k-1)m^{3}+3m^{2}+3m, & \text{if \ensuremath{p=3};}\\
3(k-1)m^{3}+3m, & \text{otherwise,}
\end{cases}
\]
and $N'=3(k-1)m^{3}$.
\item If $k$ is odd, then
\begin{align*}
N & =\begin{cases}
(3k-2)m^{3}+6m^{2}+3m, & \text{if }p=2\text{ and }\ensuremath{3\mid k};\\
(3k-2)m^{3}+3m, & \text{if }p=2\text{ and }\ensuremath{3\nmid k};\\
(2k-1)m^{3}+3m^{2}+2m, & \text{if \ensuremath{2\in\Phi};}\\
(2k-1)m^{3}+2m, & \text{otherwise,}
\end{cases}\\
\noal{and}N' & =\begin{cases}
(3k-2)m^{3} & \text{ if }p=2;\\
(2k-1)m^{3} & \text{ otherwise.}
\end{cases}
\end{align*}
\end{enumerate}
\end{thm}

A proof of the theorem will be given in the last section after we
have developed the necessary materials in the next two sections. Although
the only new part of the theorem concerns the odd case, we provide
in this paper a unified view using block intersections for both even
and odd cases.
\begin{rem}
There are prime power $q$ such that $2k\mid(q-1)$ with $(q,k)$
circular but not $(q,2k)$. For example, when $k=3$, every prime
power $q=p^{r}$ ($p$ a prime and $r\geq1$) with $3\mid(q-1)$ gives
a circular pair $(q,3)$. Now, if $p\in\{2,3,7,13,19\}$ and $6$
also divides $q-1$, then $(q,6)$ is not circular (cf.~\cite{Modisett89}
as well as \cite[Section~5]{Clay92}).
\end{rem}

\begin{rem}
The condition that $(p,k)$ be circular strange as it may seem
is actually a common phenomenon in finite fields. As shown by Modisett~\cite{Modisett89},
for each $k$, there is a finite set ${\cal P}_{k}$ of primes such
that a pair $(p,k)$ is circular if and only if $p$ is not in ${\cal P}_{k}$.
For example, ${\cal P}_{4}=\{2,5\}$, so for any prime $p$ where
$p\not=2$ and $p\not=5$, $(p,4)$ is a circular pair. Also, the
set $\{2,3,7,13,19\}$ in the previous remark is actually ${\cal P}_{6}$.
Available data shows that if we compare the size $n_{k}=|{\cal P}_{k}|$
of ${\cal P}_{k}$ to the number $\pi_{k}$ of primes less than or
equal to $p_{k}=\max{\cal P}_{k}$, then $n_{k}/\pi_{k}$ is usually
very small (see Table~\ref{maxp}).
\begin{table}
\begin{centering}
{\scriptsize{}{}{}}{\scriptsize{}%
\begin{tabular}{|c|c|c|c|c||c|c|c|c|c|}
\hline
 & {\scriptsize$p_{k}$} & {\scriptsize$\pi_{k}$} & {\scriptsize$n_{k}$} & {\scriptsize$n_{k}/\pi_{k}$\vphantom{$\Big|$}} &  & {\scriptsize$p_{k}$} & {\scriptsize$\pi_{k}$} & {\scriptsize$n_{k}$} & {\scriptsize$n_{k}/\pi_{k}$}\tabularnewline
\hline
{\scriptsize$k=4$} & {\scriptsize$5$} & {\scriptsize$3$} & {\scriptsize$2$} & {\scriptsize$0.666667$\vphantom{$\Big|$}} & {\scriptsize$k=19$} & {\scriptsize$174763$} & {\scriptsize$15898$} & {\scriptsize$54$} & {\scriptsize$0.003397$\vphantom{$\Big|$}}\tabularnewline
{\scriptsize$k=5$} & {\scriptsize$11$} & {\scriptsize$5$} & {\scriptsize$2$} & {\scriptsize$0.4$\vphantom{$\Big|$}} & {\scriptsize$k=20$} & {\scriptsize$23321$} & {\scriptsize$2600$} & {\scriptsize$58$} & {\scriptsize$0.022308$\vphantom{$\Big|$}}\tabularnewline
{\scriptsize$k=6$} & {\scriptsize$19$} & {\scriptsize$8$} & {\scriptsize$5$} & {\scriptsize$0.625$\vphantom{$\Big|$}} & {\scriptsize$k=21$} & {\scriptsize$8171731$} & {\scriptsize$550533$} & {\scriptsize$93$} & {\scriptsize$0.000169$\vphantom{$\Big|$}}\tabularnewline
{\scriptsize$k=7$} & {\scriptsize$43$} & {\scriptsize$14$} & {\scriptsize$4$} & {\scriptsize$0.285714$\vphantom{$\Big|$}} & {\scriptsize$k=22$} & {\scriptsize$165749$} & {\scriptsize$15157$} & {\scriptsize$59$} & {\scriptsize$0.003893$\vphantom{$\Big|$}}\tabularnewline
{\scriptsize$k=8$} & {\scriptsize$41$} & {\scriptsize$13$} & {\scriptsize$5$} & {\scriptsize$0.384615$\vphantom{$\Big|$}} & {\scriptsize$k=23$} & {\scriptsize$2796203$} & {\scriptsize$203095$} & {\scriptsize$78$} & {\scriptsize$0.000384$\vphantom{$\Big|$}}\tabularnewline
{\scriptsize$k=9$} & {\scriptsize$271$} & {\scriptsize$58$} & {\scriptsize$7$} & {\scriptsize$0.12069$\vphantom{$\Big|$}} & {\scriptsize$k=24$} & {\scriptsize$28753$} & {\scriptsize$3132$} & {\scriptsize$89$} & {\scriptsize$0.028416$\vphantom{$\Big|$}}\tabularnewline
{\scriptsize$k=10$} & {\scriptsize$101$} & {\scriptsize$26$} & {\scriptsize$8$} & {\scriptsize$0.307692$\vphantom{$\Big|$}} & {\scriptsize$k=25$} & {\scriptsize$9430951$} & {\scriptsize$629307$} & {\scriptsize$123$} & {\scriptsize$0.000195$\vphantom{$\Big|$}}\tabularnewline
{\scriptsize$k=11$} & {\scriptsize$683$} & {\scriptsize$124$} & {\scriptsize$8$} & {\scriptsize$0.064516$\vphantom{$\Big|$}} & {\scriptsize$k=26$} & {\scriptsize$926537$} & {\scriptsize$73227$} & {\scriptsize$111$} & {\scriptsize$0.001516$\vphantom{$\Big|$}}\tabularnewline
{\scriptsize$k=12$} & {\scriptsize$193$} & {\scriptsize$44$} & {\scriptsize$15$} & {\scriptsize$0.340909$\vphantom{$\Big|$}} & {\scriptsize$k=27$} & {\scriptsize$34975153$} & {\scriptsize$2145358$} & {\scriptsize$185$} & {\scriptsize$0.000086$\vphantom{$\Big|$}}\tabularnewline
{\scriptsize$k=13$} & {\scriptsize$2731$} & {\scriptsize$399$} & {\scriptsize$14$} & {\scriptsize$0.035088$\vphantom{$\Big|$}} & {\scriptsize$k=28$} & {\scriptsize$2968337$} & {\scriptsize$214686$} & {\scriptsize$149$} & {\scriptsize$0.000694$\vphantom{$\Big|$}}\tabularnewline
{\scriptsize$k=14$} & {\scriptsize$1289$} & {\scriptsize$209$} & {\scriptsize$20$} & {\scriptsize$0.095694$\vphantom{$\Big|$}} & {\scriptsize$k=29$} & {\scriptsize$217108153$} & {\scriptsize$11972010$} & {\scriptsize$182$} & {\scriptsize$0.000015$\vphantom{$\Big|$}}\tabularnewline
{\scriptsize$k=15$} & {\scriptsize$46261$} & {\scriptsize$4784$} & {\scriptsize 34} & {\scriptsize$0.007107$\vphantom{$\Big|$}} & {\scriptsize$k=30$} & {\scriptsize$56941$} & {\scriptsize$5775$} & {\scriptsize$134$} & {\scriptsize$0.023203$\vphantom{$\Big|$}}\tabularnewline
{\scriptsize$k=16$} & {\scriptsize$2129$} & {\scriptsize$320$} & {\scriptsize$24$} & {\scriptsize$0.075$\vphantom{$\Big|$}} & {\scriptsize$k=31$} & {\scriptsize$1114506049$} & {\scriptsize$56359124$} & {\scriptsize$257$} & {\scriptsize$0.000005$\vphantom{$\Big|$}}\tabularnewline
{\scriptsize$k=17$} & {\scriptsize$43691$} & {\scriptsize$4552$} & {\scriptsize$34$} & {\scriptsize$0.007469$\vphantom{$\Big|$}} & {\scriptsize$k=32$} & {\scriptsize$21821249$} & {\scriptsize$1378629$} & {\scriptsize$273$} & {\scriptsize$0.000198$\vphantom{$\Big|$}}\tabularnewline
{\scriptsize$k=18$} & {\scriptsize$5779$} & {\scriptsize$758$} & {\scriptsize$31$} & {\scriptsize$0.0408973$\vphantom{$\Big|$}} &  &  &  &  & \tabularnewline
\hline
\multicolumn{10}{c}{\begin{cellvarwidth}[t]
\centering
{\scriptsize{} $p_{k}=\max{\cal P}_{k}$, $\pi_{k}={}$the number
of primes less than or equal to $p_{k}$, $n_{k}=|{\cal P}_{k}|$.\vphantom{$\Big|$}}{\scriptsize\par}

{\scriptsize (Source: \href{https://nrfort.epizy.com/wiki/MathWiki/PkMax}{https://nrfort.epizy.com/wiki/MathWiki/PkMax})}
\end{cellvarwidth}}\tabularnewline
\end{tabular}}{\tiny\smallskip{}
}{\tiny\par}
\par\end{centering}
{\tiny\caption{{\small Sizes of ${\cal P}_{k}$ and ratios}}
\label{maxp}}{\tiny\par}
\end{table}

Thus, given $k$, any randomly chosen prime $p<p_{k}$ where $p$
does not divide $k$ is very likely not in ${\cal P}_{k}$ and yields
a circular pair $(p,k)$. Of course, any prime $p$ larger than $p_{k}$
makes the pair $(p,k)$ circular.
\end{rem}

\begin{rem}
In \cite{KorchmarosS99}, Korchm\'aros and Sz\"onyi also use a geometric
approach to find the number of points on the Fermat curve given by
the equation $x^{m}+y^{m}+z^{m}=0$, which we have dealt with in \cite{KK95}
and \cite{K94a} as well. They assume that the group $\Phi$ contains
the multiplicative group of some subfield $K$ of $F$.

They also showed at the end of their paper that there are only very
few common cases between our results and the results in \cite{KorchmarosS99}.
See also Lemma~\ref{l:fieldinPhi} and the comments that follow.
\end{rem}

\section{Circular pairs}

With the above notation let $\zeta$ be a generator of the cyclic
group $F^{*}$ of nonzero elements of $F$. Also, we assume $k\geq3$
and let $\Phi$ the multiplicative subgroup of $F$ of order $k$
generated by $\varphi=\zeta^{m}$. We will use $\bk=\{1,2,\dots,k-1\}$
and $\bk_{0}=\{0,1,\dots,k-1\}$.

Basic characterizations of circularity of $(F,\Phi)$ is given in
\cite[Theorem~4]{Modisett89} (see also \cite[(5.19)]{Clay92}).
\begin{lem}
\label{l:circ}The pair $(F,\Phi)$ is circular if and only $(\varphi^{i}-1)(\varphi^{i'}-1)-(\varphi^{j}-1)(\varphi^{j'}-1)\not=0$
for all \textup{$i,j,i',j'\in\bk$} with $(i,j)\not=(i',j')$ and
$(i,i')\not=(j,j')$. This is equivalent to $(\alpha-1)(\beta-1)-(\gamma-1)(\delta-1)\not=0$
for all $\alpha,\beta,\gamma,\delta\in\Phi\setminus\{1\}$ with $(\alpha,\beta)\not=(\gamma,\delta)$
and $(\alpha,\gamma)\not=(\beta,\delta)$.
\end{lem}

The statement of Theorem~\ref{Thm1_KK95} already indicates that
the question whether $2\in\Phi$ or not, is of special interest.
\begin{lem}
\label{l:2inPhi} Let $(p,k)$ be circular with $k$ even. Then $2\in\Phi$
if and only if $p=3$.
\end{lem}

\begin{proof}
Indeed, $k$ even puts $-1$ into $\Phi$. If $2\in\Phi$, then
\[
2=1+1,~-1=-2+1\ \text{and}~2^{-1}=-2^{-1}+1
\]
are all in $\Phi\cap(\Phi+1)$. The elements $2$, $-1$ and $2^{-1}$
are all distinct unless $p=3$.

The converse is clear, as $2=-1$ if the characteristic of $F$ is
$3$.
\end{proof}
The situation is quite different for odd $k$.

In general, if $p$ is a prime divisor of $2^{k}-1$ and $p$ is not
in ${\cal P}_{k}$, then $(p,k)$ is circular and $2$ is in $\Phi$.
For example, consider $k=35$. We have $2^{35}-1=31\cdot71\cdot127\cdot122921$
with $31$, $127$ and $122921$ not in ${\cal P}_{35}$. Therefore,
$(31,35)$, $(127,35)$, and $(122921,35)$ are circular pairs and
their respective $\Phi$ contains $2$. In the case of $(31,35)$,
$2$ has a multiplicative order of $5$ in every field $\text{GF}(31^{t})$
where $t$ is divisible by $6$. In the case of $(127,35)$, $2$
has a multiplicative order of $7$ in every field $\text{GF}(127^{t})$
where $t$ is divisible by $4$. In the case of $(122921,35)$, $2$
generates $\Phi$ in every field $\text{GF}(122921^{t})$ for all
$t\in\mathbb{N}$. In case when $p=2^{k}-1$ is itself a prime (i.e.,
$p$ is a Mersenne prime), $(p,k)$ is always circular (\cite{ClayYeh94},
see also \cite[Theorem 6.42]{Clay92}).

We will now establish some technical lemmas.
\begin{lem}
\label{l:p+x=00003D00003D00003D2} Let $(p,k)$ be circular. If $\lambda,\chi\in\Phi$
are such that $\lambda+\chi=2$, then $\lambda=\chi$, and in case
that $p\not=2$, one has $\lambda=\chi=1$.
\end{lem}

\begin{proof}
The case $p=2$ is trivial since $\lambda+\chi=0$ clearly implies
$\lambda=\chi$. Thus now we deal with the case $p\not=2$.

If $\chi=1$, then $\lambda=1=\chi$. So we only consider $\lambda\not=1$
and $\chi\not=1$, and assume $\lambda\not=\chi$. In such case, we
find three distinct elements $\lambda\chi^{-1}+1=2\chi^{-1}$, $\lambda^{-1}\chi+1=2\lambda^{-1}$
and $1+1=2$ in $(\Phi+1)\cap2\Phi$. This makes $(F,\Phi)$ not circular,
a contradiction.
\end{proof}
\begin{lem}
\label{l:twocase} Let $(F,\Phi)$ be circular, and let $\chi=(\psi-1)^{-1}(\lambda-1)$
where $\lambda,\psi\in\Phi\setminus\{1\}$. If $\chi\in\Phi$, then
\[
\text{ either }\chi=1\text{ and }\psi=\lambda,\text{\quad or ~}\chi=-\lambda\text{ and }\psi=\lambda^{-1}.
\]
The second case implies either $p=2$, or that $|\Phi|$ is even.
\end{lem}

\begin{proof}
From the definition of $\chi$, we have $\lambda-1-\chi(\psi-1)=0$.
If $\chi\psi=1$, then we have $\lambda+\chi=2$. As $\lambda\not=1$,
Lemma~\ref{l:p+x=00003D00003D00003D2} implies $p=2$ and $\chi=\lambda=-\lambda$.
Also, if $\chi\psi=-1$, then $\chi=-\lambda$.

Now we can assume $(\chi\psi)^{2}\not=1$. We compute
\begin{align*}
0 & =\bigl(\lambda-1-\chi(\psi-1)\bigr)(\chi^{-1}\psi^{-1}-1)\\
 & =(\lambda-1)(\chi^{-1}\psi^{-1}-1)-1+\psi^{-1}+\chi\psi-\chi\\
 & =(\lambda-1)(\chi^{-1}\psi^{-1}-1)-(\chi\psi-1)(\psi^{-1}-1).
\end{align*}
As $(F,\Phi)$ is circular the contrapostion of Lemma~\ref{l:circ}
enforces one of the following cases
\begin{enumerate}
\item[(i)] $\lambda=\psi^{-1}$ leads to $\chi=(\lambda-1)(\lambda^{-1}-1)^{-1}=-\lambda$;
\item[(ii)] $\lambda=\chi\psi$, as well as $\chi^{-1}\psi^{-1}=\psi^{-1}$,
implies $\chi=1$;
\item[(iii)] $\chi^{-1}\psi^{-1}=\chi\psi$ yields $(\chi\psi)^{2}=1$, contradicting
our assumption.
\end{enumerate}
{}Finally, if $p\not=2$, then $-\lambda=\chi\in\Phi$ puts $-1$
into $\Phi$ as $\lambda\in\Phi$. Therefore, $\Phi$ must have an
even number of elements in this case.
\end{proof}
We can generalize this lemma, and put it into a more convenient form
for later use. To do this, we define $\cC=\big\{ c_{j,i}\mid i,j\in\bk,i\not=j\big\}$
with $c_{j,i}=(\varphi^{j}-1)^{-1}(\varphi^{i}-1)$. We will use a
convention for the indices to $c_{j,i}$. Since they occur as exponents
of $\varphi$, they can be changed modulo $k$. Also, notice that
$i':=j-i$ if and only if $i=j-i'$. Therefore, we can always assume
the representative is chosen from $\bk_{0}$, even if $j-i'$ is not
in this set.
\begin{lem}
\label{l:cjiPhi} Let $(p,k)$ be circular. For $i,i',j\in\bk$ we
have
\begin{enumerate}
\item If $k$ is even, then $c_{j,i}\in\Phi c_{j,i'}\iff i=i'\mbox{ or }i=k-i'$.
\itemsep=3pt
\item If $k$ is odd, then $c_{j,i}\in\Phi c_{j,i'}\iff i=i'$ or, in case
$p=2$, $i=k-i'$.
\end{enumerate}
\end{lem}

\begin{proof}
For $\chi\in\Phi$, we have
\begin{align*}
c_{j,i}=\chi c_{j,i'} & \iff(\vp^{j}-1)^{-1}(\vp^{i}-1)=\chi(\vp^{j}-1)^{-1}(\vp^{i'}-1)\\
 & \iff(\vp^{i}-1)=\chi(\vp^{i'}-1)\\
 & \iff\frac{\vp^{i}-1}{\vp^{i'}-1}=\chi\in\Phi.
\end{align*}
The result then follows from Lemma~\ref{l:twocase}.
\end{proof}
Finally, we give a criterion for (non-)circularity. The basic idea
to the following lemma is in~\cite{KorchmarosS99}.
\begin{lem}
\label{l:fieldinPhi} Let $\Phi$ contain the multiplicative group
of some subfield $K$ of $F$. If $(F,\Phi)$ is circular, then $|K|\le4$.
\end{lem}

Note that the condition of the lemma means that $(p^{h}-1)\mid k$
for some $h\in\mathbb{N}$.
\begin{proof}
As $K^{*}\cap(K^{*}+1)=K\setminus\{0,1\}$ is contained in the intersection
$\Phi\cap(\Phi+1)$, if $(F,\Phi)$ is circular, then $|K\setminus\{0,1\}|\leq2$,
and so $|K|\leq4$.
\end{proof}
The above lemma states that when $\Phi$ contains the multiplicative
group of a subfield $K$ of $F$, the characteristic of $F$ can only
be $2$ or $3$. If $k$ is even, then $F$ has characteristic $3$.
There are indeed circular pairs $(3,k)$, such as $(3,4)$ and $(3,10)$.
However, if $k$ is divisible by $8$, then the multiplicative group
of the field with $9$ elements is contained within $\Phi$, and so
$(3,k)$ cannot be circular. This indicates that there is only a small
overlap between our results and those in~\cite{KorchmarosS99}, as
mentioned at the end of the Introduction.

\section{Block intersections}

The solutions of the equation $x^{m}+y^{m}-z^{m}=1$ over the Galois
field $F=\text{GF}(q)$ can be related to certain block intersections
$B\cap C$ for $B,C\in{\cal B}=\{\Phi a+c\mid a\in F^{*},c\in F\}$.
\begin{rem}
The pair $(F,{\cal B})$ turns out to be a $2$-design (balanced incomplete
block design) with parameters $(q,mq,q-1,k,k-1)$. If the intersection
of any two distinct blocks $B$ and $C$ in ${\cal B}$ has at most
$2$ elements, the $2$-design $(F,{\cal B})$ is said to be \textit{circular}
as the blocks behave like circles. See \cite{Clay88,Clay92}. Note
that circular $2$-designs are also referred to as ``super-simple''
in other contexts (see \cite{CharlesJ} for references).
\end{rem}

Firstly, we deal with some special intersections that will come in
handy later. For $a,c\in F^{*}$, let us denote
\[
[a,c]_{k}=|\Phi\cap(\Phi a+c)|.
\]

\begin{lem}
\label{l:propBrack} Let $a,c\in F^{*}$. For all $\psi,\chi\in\Phi$
we have
\[
[a,c]_{k}=[\psi a,\chi c]_{k}=[c,a]_{k}=[-a^{-1}c,a^{-1}]_{k}.
\]
\end{lem}

\begin{proof}
The first equality is clear. For $\lambda,\mu\in\Phi$, we have
\[
\lambda=\mu a+c\iff\lambda\mu^{-1}=\mu^{-1}c+a\iff\lambda^{-1}\mu=\lambda^{-1}(-a^{-1}c)+a^{-1},
\]
which yield the other equalities.
\end{proof}
\begin{thm}
\label{th:cycloA} Let $a,c\in F^{*}$. We have
\begin{enumerate}
\item $[a,c]_{k}>0\iff a\in\Phi(\psi-c)$ for some $\psi\in\Phi\setminus\{c\}$;\itemsep=3pt
\item $[a,c]_{k}\ge2\iff a\in\Phi c_{j,i}$ and $c\in\Phi c_{j,j-i}$ for
some $i,j\in\bk$, $i\not=j$.
\end{enumerate}
\end{thm}

\begin{proof}
The first statement is obvious.

If $[a,c]_{k}\ge2$, then there exist $\lambda_{1},\lambda_{2},\mu_{1},\mu_{2}\in\Phi$,
$\lambda_{1}\not=\lambda_{2}$ with $\lambda_{1}=\mu_{1}a+c$ and
$\lambda_{2}=\mu_{2}a+c$. Then $\mu_{1}\not=\mu_{2}$, and
\begin{alignat*}{1}
a & =(\lambda_{2}-\lambda_{1})(\mu_{2}-\mu_{1})^{-1}=\lambda_{1}\mu_{1}^{-1}(\lambda_{2}\lambda_{1}^{-1}-1)(\mu_{1}^{-1}\mu_{2}-1)^{-1}\\
\noal{and}c & =\lambda_{1}-\mu_{1}a=(\lambda_{1}\mu_{2}-\lambda_{2}\mu_{1})(\mu_{2}-\mu_{1})^{-1}=\lambda_{2}(\lambda_{1}\lambda_{2}^{-1}\mu_{1}^{-1}\mu_{2}-1)(\mu_{1}^{-1}\mu_{2}-1)^{-1}.
\end{alignat*}
Let $\mu_{1}^{-1}\mu_{2}=\varphi^{j_{0}}$ and $\lambda_{2}\lambda_{1}^{-1}=\varphi^{i_{0}}$,
where $i_{0},j_{0}\in\bk$. Then $a\in\Phi c_{j_{0},i_{0}}$ and $c\in\Phi c_{j_{0},j_{0}-i_{0}}$.
Since $c\not=0$, we have $i_{0}\not=j_{0}$.

For the converse, we first take arbitrary distinct $i,j\in\bk$ and
claim that $[c_{j,i},c_{j,j-i}]\geq2$. Direct computations give
\begin{align*}
\varphi^{-i}c_{j,i}+c_{j,j-i} & =\left(\varphi^{-i}(\varphi^{i}-1)+\varphi^{j-i}-1\right)(\varphi^{j}-1)^{-1}=(\varphi^{j-i}-\varphi^{-i})(\varphi^{j}-1)^{-1}=\varphi^{-i}\\
\noal{and}\varphi^{j-i}c_{j,i}+c_{j,j-i} & =\left(\varphi^{j-i}(\varphi^{i}-1)+\varphi^{j-i}-1\right)(\varphi^{j}-1)^{-1}=(\varphi^{j}-1)(\varphi^{j}-1)^{-1}=1.
\end{align*}
Hence
\[
\{1,\varphi^{-i}\}=\{\varphi^{j-i}c_{j,i}+c_{j,j-i},\varphi^{-i}c_{j,i}+c_{j,j-i}\}\subseteq\Phi\cap(\Phi c_{j,i}+c_{j,j-i}).
\]
This establishes the claim. Now, the first equality of Lemma~\ref{l:propBrack}
is exactly what we need for the converse. We have~(2).
\end{proof}
From Theorem~\ref{th:cycloA} we directly obtain
\begin{thm}
\label{th:cycloC} Let $(p,k)$ be circular and $a,c\in F$ with $a\not=0$.
Then
\[
[a,c]_{k}=2\iff a\in\Phi c_{j,i}\text{ and }c\in\Phi c_{j,j-i}\text{ for some }i,j\in\bk,i\not=j.
\]
\end{thm}

We look closely to the case when $c=1$, which will be needed in the
next section.
\begin{thm}
\label{th:cyclo1} Let $(p,k)$ be circular and $a\in F^{*}$. If
$k$ is even then
\begin{align*}
[a,1]_{k} & =\begin{cases}
2, & \mbox{ for }a\in\Phi(\psi-1),\psi\in\Phi\setminus\{1,-1\};\\
1, & \mbox{ for }a\in2\Phi;\\
0, & \mbox{ otherwise. }
\end{cases}\\
\noal{\text{If }k\text{ is odd, then}}[a,1]_{k} & =\begin{cases}
2, & \mbox{ for }p=2\text{ and }a\in\Phi(\psi-1),\psi\in\Phi\setminus\{1\};\\
1, & \mbox{ for }p\not=2\text{ and }a\in\Phi(\psi-1),\psi\in\Phi\setminus\{1\};\\
0, & \mbox{ otherwise. }
\end{cases}
\end{align*}
\end{thm}

\begin{proof}
By Theorem~\ref{th:cycloC}, we have $[a,1]_{k}=2$ if and only if
$a\in\Phi c_{j,i}$ and $1\in\Phi c_{j,j-i}$ for some $i,j\in\bk$,
$i\not=j$.

The condition $1\in\Phi c_{j,j-i}$ is equivalent to $c_{j,j-i}\in\Phi$,
and by Lemma~\ref{l:cjiPhi}, there are two cases. Firstly, $j=j-i$,
contradicting $i\not=0$. Secondly, $j-i=k-j$ and $2\mid k$ or $p=2$.

If $k$ is even or $p=2$, we get that $[a,1]_{k}=2$ if and only
if $i=2j$, $j\not=k/2$, and $a\in\Phi c_{j,2j}$. In this case,
$-1\in\Phi$, and we can compute
\[
a\in\Phi c_{j,2j}=\Phi(\varphi^{j}+1)=\Phi(-1)(-\varphi^{j}-1)=\Phi(\psi-1),
\]
where $\psi=-\varphi^{j}\in\Phi\setminus\{1,-1\}$. This accounts
for the two element cases.

Theorem~\ref{th:cycloA}(1) yields $[a,1]_{k}\ge1$ if and only if
$a\in\Phi(\psi-1)$ for some $\psi\in\Phi\setminus\{1\}$. This settles
the odd case with $p\not=2$.

To obtain $[a,1]_{k}=1$ in the even case we need to make sure that
$[a,1]_{k}\not=2$. Thus we are left with $a\in\Phi(\psi-1)$ and
$\psi\in\{1,-1\}$. As $\psi=1$ is not possible, the theorem is proved.
\end{proof}
\begin{rem}
Recall that for $h,\ell\in\mathbb{Z}$, the \emph{cyclomic number}
$(h,\ell)_{m}$ is the number of pairs of integers $(u,v)$, $0\le u<k$
and $0\leq v<k$, which solve the equation
\begin{equation}
\zeta^{mu+h}+1=\zeta^{mv+\ell}.\label{eq:cyclo}
\end{equation}
By setting $a=\zeta^{\ell}$, $b=\zeta^{h}$, $x=\zeta^{v}$ and $y=\zeta^{u}$,
we see that finding a solution $(u,v)$ of~(\ref{eq:cyclo}) is equivalent
to finding solutions $(x,y)$ over $F$ for the equation
\begin{equation}
ax^{m}=by^{m}+1.\label{eq:xmym}
\end{equation}
Notice that each solution $(u,v)$ of~(\ref{eq:cyclo}) gives $m^{2}$
solutions $(x,y)$ of~(\ref{eq:xmym}), namely all the elements of
the set $\{(\zeta^{v+im},\zeta^{u+jm})\mid0\leq i<m,0\leq j<m\}$.
Rewriting~(\ref{eq:xmym}) into
\[
x^{m}=a^{-1}by^{m}+a^{-1}
\]
leads to
\[
(h,\ell)_{m}=[a^{-1}b,a^{-1}]_{k}=[-b,a]_{k}
\]
by Lemma~\ref{l:propBrack}.

For fixed $c\in F$ it is clear that $\{\Phi r+c\mid r\in F^{*}\}$
is a partition of $F\setminus\{c\}$. Let $R$ be a set of representatives
of the cosets of $\Phi$. Thus, from Lemma~\ref{l:propBrack} one
derives again the following identities (in \cite[eq.~(17), p.~394]{Dickson 35a}):
\begin{align*}
\sum_{r\in R}[r,c]_{k} & =\begin{cases}
k-1, & \mbox{ for }c\in\Phi;\\
k, & \mbox{ for }c\not\in\Phi,
\end{cases}\\
\noal{and}\sum_{r\in F^{*}}[r,c]_{k} & =\begin{cases}
k(k-1), & \mbox{ for }c\in\Phi;\\
k^{2}, & \mbox{ for }c\not\in\Phi.
\end{cases}
\end{align*}
\end{rem}

Now we go directly for the main theorem. Let us put
\begin{align*}
S & =\{(x,y,z)\in F^{3}\mid x^{m}+y^{m}-z^{m}=1\text{ and }xyz\not=0\},\\
\noal{and}T & =\{(x,y,z)\in F^{3}\mid x^{m}+y^{m}-z^{m}=1,\;xyz=0\}.
\end{align*}
Thus, the disjoint union $S\cup T$ is the complete set of the solutions
of the equation (\ref{eq:maxp}) in $F^{3}$, while $S$ is the set
of those in~$F^{*}{}^{3}$.

Consider $S$ first. Suppose that there exists $w\in(\Phi+1)\cap(\Phi+\varphi^{i})$
for some $i\in\bk_{0}$. Then there are $u,v\in\bk_{0}$ such that
$\varphi^{v}+1=w=\varphi^{u}+\varphi^{i}$. Recall that $\zeta$ is
a generator $F^{*}$ with $\varphi=\zeta^{m}$. For arbitrary $m$th
roots of unity $\lambda_{1},\lambda_{2},\lambda_{3}\in\langle\zeta^{k}\rangle$,
put $x=\lambda_{1}\zeta^{i}$, $y=\lambda_{2}\zeta^{u}$ and $z=\lambda_{3}\zeta^{v}$,
and we have
\[
x^{m}+y^{m}-z^{m}=\lambda_{1}^{m}\zeta^{m\cdot i}+\lambda_{2}^{m}\zeta^{m\cdot u}-\lambda_{3}^{m}\zeta^{m\cdot v}=\varphi^{i}+\varphi^{u}-\varphi^{v}=1.
\]
Therefore $(x,y,z)\in S$.

Conversely, let $(x,y,z)\in S$. Then $x=\zeta^{i'}$, $y=\zeta^{u'}$,
and $z=\zeta^{v'}$ for some $i',u',v'\in\mathbb{Z}$. There exist
$i'',u'',v''\in\mathbb{Z}$ and $i,u,v\in\bk_{0}$ with $i'=i''k+i$,
$u'=u''k+u$ and $v'=v''k+v$, and so $x^{m}=\zeta^{i\cdot m}=\varphi^{i}$,
$y^{m}=\zeta^{u\cdot m}=\varphi^{u}$, and $z^{m}=\zeta^{v\cdot m}=\varphi^{v}$.
From $x^{m}+y^{m}-z^{m}=1$, we thus get
\[
1+\varphi^{v}=\varphi^{u}+\varphi^{i}\in(\Phi+1)\cap(\Phi+\varphi^{i}).
\]

Therefore $|S|$ is given by
\begin{align}
|S| & =m^{3}\cdot\sum_{i=0}^{k-1}\big|(\Phi+1)\cap(\Phi+\varphi^{i})\big|\nonumber \\
 & =m^{3}\cdot\left(k+\sum_{i=1}^{k-1}\big|(\Phi+1)\cap(\Phi+\varphi^{i})\big|\right).\label{eq:S-1}
\end{align}

Now, consider $T$. Let
\begin{alignat*}{2}
T_{x} & =\{(0,y,z)\in T\mid yz\not=0\},\text{ } & T_{y} & =\{(x,0,z)\in T\mid xz\not=0\},\\
T_{z} & =\{(x,y,0)\in T\mid xy\not=0\},\\
T_{x,y} & =\{(x,y,z)\in T\mid x=y=0\},\text{ } & T_{x,z} & =\{(x,y,z)\in T\mid x=z=0\},\\
T_{y,z} & =\{(x,y,z)\in T\mid y=z=0\}.
\end{alignat*}
Then $T$ is the disjoint union of $T_{x}$, $T_{y}$, $T_{z}$, $T_{x,y}$,
$T_{x,z}$, and $T_{y,z}$. Using similar arguments to that for $S$,
we can see that
\[
|T_{x}|=|T_{y}|=m^{2}\cdot|\Phi\cap(\Phi+1)|\text{ and }|T_{z}|=m^{2}\cdot|\Phi\cap(1-\Phi)|.
\]
From Lemma~\ref{l:propBrack} we find $|\Phi\cap(1-\Phi)|=[-1,1]_{k}=[1,1]_{k}=|\Phi\cap(\Phi+1)|$.

Finally, it is clear that
\[
|T_{yz}|=|\{x\in F\mid x^{m}=1\}|=m\text{ and }|T_{xz}|=|\{y\in F\mid y^{m}=1\}|=m.
\]
Since the number of the solutions of the equation $-z^{m}=1$ in $F$
is $m$ if $-1\in\Phi$ and $0$ otherwise, $|T_{xy}|$ is $m$ or
$0$ depending on whether $-1$ is in $\Phi$ or not. Putting the
above together, and note that $-1\in\Phi$ if and only if $|\Phi|$
is even or $p=2$, we have
\begin{equation}
|T|=\begin{cases}
3m^{3}\cdot|\Phi\cap(\Phi+1)|+3m, & \text{if }k\text{ is even or }p=2;\\
3m^{3}\cdot|\Phi\cap(\Phi+1)|+2m, & \text{if }k\text{ is odd and }p\not=2.
\end{cases}\label{eq:T-1}
\end{equation}

Now, set
\begin{equation}
s=\sum_{j=1}^{k-1}\big|(\Phi+1)\cap(\Phi+\varphi^{j})\big|\label{eq:s}
\end{equation}
and
\begin{equation}
t=[1,1]_{k}=|\Phi\cap(\Phi+1)|.\label{eq:t}
\end{equation}
With (\ref{eq:S-1}) and (\ref{eq:T-1}), we have the following general
theorem.
\begin{thm}
\label{thm:general}Let $F=\textup{GF}(q)$, $k\geq3$ a divisor of
$q-1$ and $m=(q-1)/k$. Assume that $(p,k)$ is circular, then the
number of solutions of $x^{m}+y^{m}-z^{m}=1$ in $F$ is given by
\[
\begin{cases}
m^{3}(k+s)+3m^{2}t+2m, & \text{if }\ensuremath{k}\text{ is odd};\\
m^{3}(k+s)+3m^{2}t+3m, & \text{if }\ensuremath{k}\text{ is even or }p=2.
\end{cases}
\]
Here, $s$ and $t$ are as in $(\ref{eq:s})$ and $(\ref{eq:t})$.
\end{thm}

It is clear from this theorem that to prove Theorem~\ref{Thm1_KK95},
we just have to find formulas for $s$ and $t$ under the condition
that $(F,\Phi)$ is circular.

\section{Proof of Theorem \ref{Thm1_KK95}}

The following improvement of \cite[(3.3)]{KK95} gives us $t$.
\begin{lem}
\label{lem:t}Let $(p,k)$ be circular. Then
\[
(1)\quad2\mid k\Rightarrow t=\begin{cases}
2, & \text{if }6\mid k;\\
1, & \mbox{if }p=3;\\
0, & \mbox{otherwise. }
\end{cases}\hspace*{0.85cm}\hspace*{\fill}(2)\quad2\nmid k\Rightarrow t=\begin{cases}
2, & \text{when }p=2\text{ and }3\mid k;\\
1, & \mbox{when }p\not=2\text{ and }2\in\Phi;\\
0, & \mbox{otherwise. }
\end{cases}\hspace*{\fill}\hspace*{1.3cm}
\]
\end{lem}

\begin{proof}
The case when $k$ is even is already given in \cite[(3.3)(1)]{KK95}.

Let $k$ be odd. Consider $p\not=2$ first. By Theorem~\ref{th:cyclo1},
we must have $t=[1,1]_{k}\leq1$. When $2\in\Phi$, we have $2=1+1\in\Phi\cap(\Phi+1)$,
and so $t=1$. Conversely, assume $t=1$ and let $\psi,\chi\in\Phi$
such that $\psi=\chi+1\in\Phi\cap(\Phi+1)$. From this, we have $\psi\chi^{-1}=\chi^{-1}+1\in\Phi\cap(\Phi+1)$,
and so $\psi=\psi\chi^{-1}$. Then $\chi=1$ and $\psi=2\in\Phi$.

Let $p=2$. We note that $3\mid k$ implies the existence of a third
root of unity $\tau\in\Phi$. Then we have $\tau^{2}=\tau+1\in\Phi\cap(\Phi+1)$
as well as $\tau=\tau^{2}+1\in\Phi\cap(\Phi+1)$. Thus, $t=2$. Conversely,
assume $t\not=0$. So $\lambda=\delta+1\in\Phi\cap(\Phi+1)$ for some
$\lambda,\delta\in\Phi\setminus\{1\}$. Then we also have
\[
\delta=\lambda+1,\delta^{-1}\lambda=\delta^{-1}+1,\delta^{-1}=\delta^{-1}\lambda+1\in\Phi\cap(\Phi+1).
\]
By circularity we must have $\lambda=\delta^{-1}$ as $\lambda=\delta$
and $\lambda=\delta^{-1}\lambda$ are both not possible. This yields
$\delta^{2}+\delta+1=0$. Therefore $\delta$ is a third root of unity,
and so $3\mid k$. Then again, $t=2$.
\end{proof}
The next result helps us find $s$.
\begin{lem}
\label{l:intersect}Let $(F,\Phi)$ be circular and $j\in\bk$.
\begin{enumerate}
\item If $k$ is even, then
\[
\big|(\Phi+1)\cap(\Phi+\varphi^{j})\big|=\begin{cases}
2, & \mbox{ if }1\in\Phi c_{j,i}\text{ for some }i\in\bk\setminus\left\{ \frac{k}{2}\right\} ;\\
1, & \mbox{ if }1\in\Phi c_{j,\frac{k}{2}};\\
0, & \mbox{ otherwise. }
\end{cases}
\]
\item If $k$ is odd, then \itemsep=3pt
\[
\hspace*{0.65cm}\big|(\Phi+1)\cap(\Phi+\varphi^{j})\big|=\begin{cases}
2, & \mbox{ if }p=2\text{ and }1\in\Phi c_{j,i}\text{ for some }i\in\bk;\\
1, & \mbox{ if }p\not=2\text{ and }1\in\Phi c_{j,i}\text{ for some }i\in\bk;\\
0, & \mbox{ otherwise. }
\end{cases}
\]
\end{enumerate}
\end{lem}

\begin{proof}
For $\chi,\psi\in\Phi$, we have $\chi+1=\psi+\varphi^{j}\in(\Phi+1)\cap(\Phi+\varphi^{j})$
if and only if
\[
1=\frac{\chi-\psi}{\varphi^{j}-1}=\psi\frac{\chi\psi^{-1}-1}{\varphi^{j}-1}\in\Phi c_{j,i},\mbox{ where }\chi\psi^{-1}=\varphi^{i}\text{ for some }i\in\bk.
\]
Notice that $j\not=0$ ensures $i\not=0$. This already accounts for
the ``otherwise''-part in both cases. So we assume that $1\in\Phi c_{j,i}$.
Then $(\Phi+1)\cap(\Phi+\varphi^{j})$ is nonempty, and
\[
\chi+1=\psi+\varphi^{j}
\]
for some $\chi,\psi\in\Phi$ with $\chi\psi^{-1}=\varphi^{i}$.

Suppose $k$ is odd. From $|(\Phi+1)\cap(\Phi+\varphi^{j})|=|\Phi\cap(\Phi+\varphi^{j}-1)|=[1,c]_{k}$,
where $c=\varphi^{j}-1\not=0$, we have
\[
\big|(\Phi+1)\cap(\Phi+\varphi^{j})\big|=[1,c]_{k}=[c,1]_{k}
\]
by Lemma~\ref{l:propBrack} and Theorem~\ref{th:cyclo1}. This finishes
the odd case.

Now, assume that $k$ is even. Then $-1\in\Phi$, and so
\[
-\psi+1=-\chi+\varphi^{j}\in(\Phi+1)\cap(\Phi+\varphi^{j}).
\]
Hence, except in the case $\psi=-\chi$, we have already two points
of intersection.

To finish the proof, we let $\psi=-\chi$ (so that $\chi\psi^{-1}=-1=\varphi^{k/2}$),
and claim that $\big|(\Phi+1)\cap(\Phi+\varphi^{j})\big|=1$. Suppose
that this is not the case. Then there is another point $\lambda+1=\mu+\varphi^{j}$
in the intersection, where $\lambda,\mu\in\Phi$ with $\lambda\not=\chi$
(hence $\mu\not=\psi$). Again, we have
\[
-\mu+1=-\lambda+\varphi^{j}\in(\Phi+1)\cap(\Phi+\varphi^{j}).
\]
From circularity, this has to be one of the two points in the intersection.
Thus, either $\mu=-\lambda$ or $\mu=-\chi=\psi$. The former case
gives $1=\psi c_{j,\frac{k}{2}}=\mu c_{j,\frac{k}{2}}$, which still
leads to $\psi=\mu$. As this is not possible, we have $\big|(\Phi+1)\cap(\Phi+\varphi^{j})\big|=1$
as claimed, and the proof is complete.
\end{proof}
\begin{lem}
\label{t:regular}Suppose that $(F,\Phi)$ is circular. Then
\[
s=\begin{cases}
2k-3, & \text{if }\ensuremath{k}\text{ is even};\\
2k-2, & \text{if }p=2;\\
k-1, & \text{if }p\not=2\text{ and }\ensuremath{k}\text{ is odd};
\end{cases}
\]
\end{lem}

\begin{proof}
We have $1\in\Phi c_{j,i}$ if and only if $c_{j,i}\in\Phi$.

When $k$ is even or $p=2$, Lemma~\ref{l:cjiPhi} gives $c_{j,i}\in\Phi$
exactly when $j=i$ or $j=k-i$. By Lemma~\ref{l:intersect} every
$j\in\bk$ other than $j=k/2$ contributes $2$ to the sum, whereas,
$j=k/2$ contributes only one. Summarizing, we get $s=2(k-2)+1=2k-3$
when $k$ is even, and $s=2(k-1)=2k-2$ if $p=2$.

When $p\not=2$ and $k$ is odd, $c_{j,i}\in\Phi$ is only possible
when $i=j$ by Lemma~\ref{l:cjiPhi}. Therefore, Lemma~\ref{l:intersect}
shows that every $j\in\bk$ contributes $1$ to the sum. Hence $s=k-1$.
\end{proof}
Now, Theorem~\ref{Thm1_KK95} follows from Theorem~\ref{thm:general}
together with Lemmas~\ref{lem:t} and~\ref{t:regular}.

\end{document}